\documentclass{ifacconf}

\usepackage[round]{natbib} 

\usepackage{amsmath,amssymb}
\usepackage{amsfonts}
\newtheorem{theorem}{Theorem}[section]

\newtheorem{lemma}[theorem]{Lemma}

\newtheorem{remark}[theorem]{Remark}

\newcommand{\rd}{{\mathbb{R}^d}}

\newcommand{\rr}{{\mathbb R}}

\def\N{{\mathbb N}}

\def\E{{\mathbb E}}

\def\D{{\mathbb D}}

\begin{document}
\begin{frontmatter}

\title{Fractional Cauchy problems on bounded domains:  survey of recent results}


\author[First]{Erkan Nane}

\address[First]{Department of Mathematics and Statistics, 221 Parker Hall, Auburn University, Auburn, Alabama 36849, USA (e-mail: ezn0001@ auburn.edu)}

\begin{abstract}
In a fractional Cauchy problem, the usual first order time
derivative is replaced by a fractional derivative.
This problem was first considered by  \citet{nigmatullin}, and  \citet{zaslavsky} in  $\mathbb R^d$  for modeling some physical phenomena.
 The fractional
derivative models time delays in a diffusion process.
We will give a  survey of the recent results on the fractional Cauchy problem and its generalizations on bounded domains $D\subset \rd$ obtained in \citet{m-n-v-aop, mnv-2}. We also study the solutions of fractional Cauchy problem where the first time derivative is replaced with an infinite sum of fractional derivatives. We point out connections to eigenvalue problems for the fractional time operators considered.  The solutions to the eigenvalue problems are  expressed by Mittag-Leffler functions and its generalized versions. The stochastic solution of the eigenvalue problems for the  fractional derivatives are given by inverse  subordinators.
\end{abstract}


\end{frontmatter}


\section{Introduction}

A celebrated paper of \citet{Einstein1906} established a mathematical link
between random walks, the diffusion equation, and Brownian motion.  The scaling limits of a simple random walk
with mean zero, finite variance jumps yields a Brownian motion.  The probability densities of the Brownian motion
 variables solve a diffusion equation, and hence we refer to the Brownian motion as the stochastic solution to the
 diffusion (heat) equation.   The diffusion equation is the most familiar
 Cauchy problem.  The general abstract Cauchy problem is $\partial_t u=L u$, where $u(t)$ takes values in a
 Banach space and $L$ is the generator of a continuous semigroup on that space, see \citet{ABHN}.  If $L$ generates a Markov process, then we call this Markov process a stochastic solution to the Cauchy problem $\partial_t u=L u$, since its probability densities (or distributions) solve the Cauchy problem.  This point of view has proven useful, for instance, in the modern theory of fractional calculus, since fractional derivatives are generators of certain ($\alpha$-stable) stochastic processes, see \citet{limitCTRW}.

Fractional derivatives are almost as old as their more familiar integer-order counterparts, see
\citet{MillerRoss, Samko}.  Fractional diffusion equations have recently been applied to problems
in physics, finance, hydrology, and many other areas, see \citet{GorenfloSurvey,koch2,MetzlerKlafter,scalas1}.
Fractional space derivatives are used to model anomalous diffusion or dispersion, where a particle plume spreads
at a rate inconsistent with the classical Brownian motion model, and the plume may be asymmetric. When a fractional
derivative replaces the second derivative in a diffusion or dispersion model, it leads to enhanced diffusion
(also called superdiffusion).  Fractional time derivatives are connected with anomalous subdiffusion, where a cloud of particles spreads more slowly than a classical diffusion.  Fractional Cauchy problems replace the integer time derivative by its fractional counterpart:  $\partial^\beta_t u=L u$.  Here, $\partial^{\beta}_t g(t)$ indicates the Caputo fractional derivative in time, the inverse Laplace transform of $s^{\beta}\tilde{g}(s)-s^{\beta -1}g(0)$, where
$\tilde{g}(s)=\int_{0}^{\infty}e^{-st}g(t)dt$ is the usual Laplace transform, see \cite{Caputo}.
\citet{nigmatullin} gave a physical derivation of the fractional Cauchy problem,
when $L$ is the generator of some continuous Markov process $\{Y(t)\}$ started at $x=0$.
The mathematical study of fractional Cauchy problems was initiated by \citet{koch1, koch2, sch-wyss}.
The existence and uniqueness of solutions was proved in \citet{koch1, koch2}.
Fractional Cauchy problems were also invented independently by
 \citet{zaslavsky} as a model for Hamiltonian chaos.

Stochastic solutions of fractional Cauchy problems are subordinated processes.  If $X(t)$ is a stochastic solution to the Cauchy problem $\partial_t u=Au$, then under certain technical conditions, the subordinate process $X(E(t))$ is a stochastic solution to the fractional Cauchy problem $\partial^\beta_t u=Au$, see \citet{fracCauchy}.  Here, $E(t)$ is the inverse or hitting time process to a stable subordinator $D(t)$ with index $\beta\in(0,1)$.  That is, $E(t)=\inf\{x>0:D(x)>t\}$, and $D(t)$ is a L\'evy process (continuous in probability with independent, stationary increments) whose smooth probability density $f_{D(1)}(t)$ has Laplace transform $e^{-s^\beta}=\tilde f_{D(1)}(s)$, see \citet{sato}.  Just as Brownian motion is a scaling limit of a simple random walk, the stochastic solution to certain fractional Cauchy problems are scaling limits of continuous time random walks, in which the independent identically distributed (iid) jumps are separated by iid waiting times, see \citet{Zsolution}.  Fractional time derivatives arise from power law waiting times, where the probability of waiting longer than time $t>0$ falls off like $t^{-\beta}$ for $t$ large, see \citet{limitCTRW}.  This is related to the fact that fractional derivatives are non-local operators defined by convolution with a power law, see \citet{fracCauchy}.

In some applications, the waiting times between particle jumps
evolve according to a more complicated process that cannot be
adequately described by a single power law.  Then, a waiting time
model that is conditionally power law leads to a distributed-order
fractional derivative in time, defined by integrating the
fractional derivative of order $\beta$ against the probability
distribution of the power-law index, see \cite{M-S-ultra}. The
resulting distributed-order fractional Cauchy problem provides a
more flexible model for anomalous sub-diffusion.  The L\'evy
measure of a stable subordinator with index $\beta$ is integrated
against the power law index distribution to define a subordinator
$W(t)$.  Its inverse $E(t)$ produces a stochastic solution $X(E(t))$
of the distributed-order fractional Cauchy problem on $\rd$, when
$X(t)$ solves the original
 Cauchy problem $\partial u/\partial t=Lu$.


\section{ Stochastic solution of heat equation on bounded domains}\label{sec2}

Let $D$ be a bounded domain in $\rd$. We denote by $C^k(D),
C^{k,\alpha}(D), C^k(\bar D)$  the space of k-times differentiable
functions in $D$, the space of $k$-times differential functions
with $k$-th derivative is H\"older continuous of index $\alpha$,
 and the space of functions that have all the derivatives up to order  $k$ extendable continuously up to the
 boundary $\partial D$ of $D$, respectively. We refer to  \citet{m-n-v-aop} for a detailed discussion of
 these spaces and concepts in this section.

Suppose that every point of $\partial D$ is regular for $D^C$.
The corresponding Markov process is a killed Brownian
motion. We denote the eigenvalues and the eigenfunctions of the Laplacian
$\Delta=\sum_{i=1}^d\partial^2_{x_i}$ by $\{\mu_n, \phi_n\}_{n=1}^\infty$, where $\phi_n\in C^{\infty}(D)$.
\begin{remark}
 In particular, the initial function $f$ regarded as an element of $L^2(D)$ can be represented as
 \begin{equation}
 f(x)=\sum_{n=1}^\infty \bar f(n)\phi_n(x);\ \ \bar f(n)=\int_D\phi_n(x)f(x)dx.
 \end{equation}

\end{remark}
The
corresponding heat kernel is given by
$$
p_D(t,x,y)=\sum_{n=1}^{\infty}e^{-\mu_n t}\phi_n(x)\phi_n(y).
$$
The series
converges absolutely and  uniformly on $[t_0,\infty)\times D\times D$ for all
$t_0>0$. In this case, the semigroup given by
\begin{equation}\label{stochastic-killed-bm}\begin{split}
T_D(t)f(x)&=E_x[f(X(t))I(
t<\tau_D(X))]\\
&=\sum_{n=1}^{\infty}e^{-\mu_n
t}\phi_n(x)\bar{f}(n)
=\int_Dp_D(t,x,y)f(y)dy
\end{split}\end{equation}
solves the heat equation in $D$ with Dirichlet boundary conditions:
\begin{eqnarray}
\partial_t u(t,x)&=& \Delta u(t,x),\ \ x\in D, \ t>0, \nonumber\\
u(t,x)&=&0,\ \ x\in \partial D,\nonumber\\
u(0,x)&=&f(x), \ \ x\in D.\nonumber
\end{eqnarray}

\section{Fractional Cauchy problem}


 Fractional derivatives in time are useful for physical models that involve sticking or trapping, see
 \citet{Zsolution}.  They are closely connected to random walk models with long waiting times between
 particle jumps, see \citet{limitCTRW}.  The fractional derivatives are essentially convolutions with a power law.
  Various forms of the fractional derivative can be defined, depending on the domain of the power law kernel,
  and the way boundary points are handled, see \citet{MillerRoss, Samko}.  The Caputo fractional derivative
   invented by \citet{Caputo} is defined for $0<\beta<1$ as
\begin{equation}\label{CaputoDef}
\partial^\beta_t u(t,x)=\frac{1}{\Gamma(1-\beta)}\int_0^t \partial_r
u(r,x)\frac{dr}{(t-r)^\beta} .
\end{equation}
Its Laplace transform
\begin{equation}\label{CaputoLT}
\int_0^\infty e^{-st} \partial^\beta_t u(t,x)\,ds=s^\beta \tilde u(s,x)-s^{\beta-1} u(0,x)
\end{equation}
incorporates the initial value in the same way as the first
derivative.  The Caputo derivative is useful for solving
differential equations that involve a fractional time
derivative, see \citet{GorenfloSurvey,Podlubny}, because it naturally incorporates initial values.

Let $D\in\rd$ be a bounded domain. In this section we will consider  the fractional Cauchy problem:
\begin{eqnarray}
\partial^\beta_t u(t,x)
&=&
\Delta u(t,x),  \  \ x\in D, \ t> 0\label{frac-heat};\\
u(t,x)&=&0, \ \ x\in \partial D \ t> 0; \nonumber\\
u(0,x)& =& f(x), \ x\in D.\nonumber
\end{eqnarray}
To obtain a solution, let $u(t,x)=G(t)F(x)$ be a solution of
(\ref{frac-heat}).
Substituting in the PDE (\ref{frac-heat}) leads
to
$$
F(x)\partial^\beta_t G(t)
= G(t)\Delta F(x)
$$
and now dividing both sides by $G(t)F(x),$ we obtain
$$
\frac{\partial^\beta_t G(t)}{G(t)} = \frac{\Delta F(x)}{F(x)}=-\mu.
$$
That is,
\begin{equation}\label{time-pde-2}
\partial^\beta_t G(t)=-\mu G(t), \ t>0;
\end{equation}
\begin{equation}\label{space-pde-2}
\Delta F(x)=-\mu F(x) \ x\in D; \ F(x)=0, \ x\in \partial D.
\end{equation}

\noindent Eigenvalue problem (\ref{space-pde-2}) is solved by an infinite
sequence of pairs $(\mu_n, \phi_n)$, $n \ge 1,$ where
 $\phi_n$ is a sequence
of functions that form a complete orthonormal set in $L^2(D)$, $\mu_1<\mu_2\leq \mu_2\leq \cdots,$ and $\mu_n\to \infty$.

\noindent Using the $\mu_n$ determined by (\ref{space-pde-2}), we need to find a solution of
 \eqref{time-pde-2} with $\mu=\mu_n$, which is the eigenvalue problem for the Caputo fractional derivative.

 We next consider the eigenvalue problem for the Caputo fractional derivative of order $0<\beta<1$.
 \begin{lemma}\label{frac-eigenvalue-lemma}
Let $\lambda>0$. The unique solution of the eigenvalue problem
\begin{equation}\label{frac-eigenvalue}
  \partial^\beta_t G(t)=-\lambda G(t), \ \ G
  (0)=1.
 \end{equation}
 is given by the Mittag-Leffler function
\begin{equation}
G(t)=M_\beta(-\lambda t^\beta)=\sum_{n=0}^\infty\frac{(-\lambda t^\beta)^n}{\Gamma(1+\beta n)}
\end{equation}
\end{lemma}

For a detailed study of the  Mittag-Leffler type functions we refer the reader to the tutorial paper  by \citet{gorenflo-mainardi}.

Therefore the  solution to \eqref{time-pde-2} is given by
$$
G(t)=G_0(n)M_\beta(-\mu t^\beta),
$$
where $G_0(n)=\bar f(n)$ is selected to satisfy the initial condition $f$.
 Therefore using this lemma, we obtain a formal solution of the fractional Cauchy problem \eqref{frac-heat} as
\begin{equation}\label{formal-sol-L-1}\begin{split}
u(t,x)&=\sum_{n=1}^{\infty}\bar{f}(n)M_\beta(-\mu_n t^\beta)\phi_n(x)
\end{split}\end{equation}

\begin{remark}
The separation of variables technique works for a large class operators including uniformly elliptic operators $L$ in divergence form, see Remark \ref{extension-L}. For details see \citet{m-n-v-aop}.
Define a  cube
$$
D=\{x=(x_1,x_2,\cdots, x_d): \ 0<x_i<M \ \mathrm{for \ all}\ 1\leq i\leq d\}.
$$
The functions
$$
\phi_n(x)=(2/M)^{d/2}\prod_{i=1}^d\sin (\pi n_i x_i/M)
$$
parametrised by the multi-index of positive integers $n=\{n_1, n_2,\cdots, n_d\}$, form a complete orthonormal set of eigenfunctions of the Laplacian,
with Dirichlet boundary conditions with corresponding eigenvalues
$$
\mu_n=\pi^2M^{-2}(n_1^2+\cdots +n_d^2).
$$
See, for example, Lemma 6.2.1 in \citet{davies-d}. In this case the boundary of the cube domain is not smooth.
\end{remark}

\subsection{Stochastic solution.}

Fractional time derivatives emerge in anomalous diffusion models,
when particles wait a long time between jumps.
 In the standard model, called a continuous time random walk (CTRW), a particle waits a random time $J_n>0$ and
  then takes a step of random size $Y_n$.  Suppose that the two sequences
  of i.i.d.\ random variables $(J_n)$ and $(Y_n)$ are independent.
  The particle
   arrives at location $X(n)=Y_1+\cdots+Y_n$ at time $T(n)=J_1+\cdots+J_n$.   Since
    $N(t)=\max\{n\geq 0:T(n)\leq t\}$ is the number of jumps by time $t>0$, the particle location at
    time $t$ is $X(N(t))$.  If $\E Y_n=0$ and $\E [Y_n^2]<\infty$ then, as the time scale $c\to\infty$, the
    random walk of particle jumps has a scaling limit $c^{-1/2}X([ct])\Rightarrow B(t)$, a standard Brownian
    motion.  If
$P(J_n>t)\sim ct^{-\beta}$ for some $0<\beta<1$ and $c>0$, then
the scaling limit $c^{-1/\beta}T([ct])\Rightarrow D(t)$ is a
strictly increasing stable L\'evy process with index $\beta$,
sometimes called a stable subordinator.  The jump times $T(n)$ and
the number of jumps $N(t)$ are inverses $\{N({t})\geq n\}=\{T(n)\leq
t\}$, and it follows that the scaling limits are also inverses, see
 \citet[Theorem 3.2]{limitCTRW}:
$c^{-\beta}N({ct})\Rightarrow E(t)$, where
\begin{equation}\label{Etdef}
E(t)=\inf\{\tau:D(\tau)> t\} ,
\end{equation}
so that $\{E(t)\leq \tau\}=\{D(\tau)\geq t\}$.  A continuous mapping argument in \citet[Theorem 4.2]{limitCTRW} yields
the CTRW scaling limit:  Heuristically, since $N(ct)\approx c^{\beta}E(t)$, we have
$c^{-\beta/2}X(N({[ct]}))\approx (c^\beta)^{-1/2}X(c^\beta
E(t))\approx B(E(t))$, a time-changed Brownian motion.  The density $u(t,x)$ of the process $B(E(t))$ solves a
fractional Cauchy problem
\[\partial^\beta_tu(t,x)=\partial^2_x u(t,x)\]
 where the order of the fractional derivative equals the index of the stable subordinator.
 Roughly speaking, if the probability of waiting longer than time $t>0$ between jumps falls off like
 $t^{-\beta}$, then the limiting particle density solves a diffusion equation that involves a fractional
  time derivative of the same order $\beta$.

The  Laplace transform of $D(t)$ is given by
$$
\E(e^{-sD(t)})=\int_0^\infty e^{-sx}f_{D(t)}(x)dx=e^{-ts^\beta}.
$$
The inverse to the stable subordinator $E(t)$ of  $D(t)$ has density
\begin{equation}\label{E-density}
\begin{split}
f_{E(t)}(l)&=\frac{\partial}{\partial l}P(E(t)\leq l)=\frac{\partial}{\partial l}(1-P(D(l)\leq t))\\
&=-\frac{\partial}{\partial l}\int_0^{\frac{t}{l^{1/\beta}}}f_{D(1)}(u)du\\
&=(t/\beta)f_{D(1)}(tl^{-1/\beta })l^{-1-1/\beta}
\end{split}\end{equation}
using the scaling property of the density $f_{D(t)}(l)=t^{-1/\beta}f_{D(1)}(lt^{-1/\beta})$, see \cite{bertoin}.

Using the representation \eqref{E-density} and taking Laplace transforms we can show that the unique solution of the eigenvalue problem \eqref{frac-eigenvalue}
 is also given by
\begin{equation}\label{stoch-rep-m}\begin{split}
G(t)
&=\int_0^\infty \exp(-l\lambda) f_{E(t)}(l)dl=\E(\exp(-\lambda E(t)))
\end{split}
\end{equation}
see \citet{Zsolution} for the details.

Now using \eqref{stochastic-killed-bm} and \eqref{stoch-rep-m}  we can express the solution to \eqref{frac-heat} as
\begin{equation}
\begin{split}
u(t,x)&=\sum_{n=1}^{\infty}\bar{f}(n)M_\beta(-\mu_n t^\beta)\phi_n(x)\\
&=\sum_{n=1}^{\infty}\bar{f}(n)[\int_0^\infty \exp(-l\mu_n) f_{E(t)}(l)dl]\phi_n(x)\\
&=\int_0^\infty[\sum_{n=1}^{\infty}\bar{f}(n)\exp(-l\mu_n)\phi_n(x)]f_{E(t)}(l)dl\\
&=\int_0^\infty [T_D(l)f(x)]f_{E(t)}(l)dl\\
&=\E_x[f(B(E(t)))I(\tau_M>E(t))]\\
&=\E_x[f(B(E(t)))I(\tau_M(X(E))>t)].
\end{split}
\end{equation}
\begin{remark}
 \citet{m-n-v-aop} established the conditions on the initial function $f$  under which $u(t,x)$ is a classical solution ( i.e., for each $t>0$, $u(t,x)\in C^1(\bar D)\cap C^2(D)$ and  for each $x\in D$, $u(t,x)\in C^1(0,\infty)$) of \eqref{frac-heat}: that $\Delta f(x)$ has an eigenfunction expansion w.r.t. $\{\phi_n\}$ that is absolutely and uniformly convergent. The analytic expression in $(0,M)\subset\mathbb R$ above is due to \citet{agrawal}.
\end{remark}

\section{Distributed-order fractional cauchy problems}\label{distributed-order-sec}

 Let $\mu$ be a finite measure with $supp \mu\subset (0,1)$. We consider the  distributed order-time fractional derivative
\begin{equation}\label{DOFDdef-D}
\D^{(\nu)}u(t,x):=\int_0^1 \partial^\beta_t u(t,x)\nu(d\beta) ,\ \nu(d\beta)=\Gamma(1-\beta)\mu(d\beta).
\end{equation}
To ensure that $\D^{(\nu)}$ is well-defined, we impose the condition
\begin{equation}\label{finite-mu-bound-D}
\int_0^1  \frac 1{1-\beta}\, \mu(d\beta)<\infty
\end{equation}
as in  \citet[Eq.\ (3.3)]{M-S-ultra}.  Since $\Gamma(x)\sim 1/x$, as $x\to 0+$, this ensures that $\nu(d\beta)$ is a finite measure on $(0,1)$.

\subsection{Eigenvalue problem: solution with waiting time process}
Stochastic solution to the distributed-order fractional Cauchy problem is
obtained by considering a a more flexible  sequence of CTRW.  At each scale
  $c>0$, we are given i.i.d.\ waiting times $(J_n^{c})$ and i.i.d.\ jumps $(Y_n^{c})$.  Assume the waiting times
   and jumps form triangular arrays whose row sums converge in distribution. Letting
   $X^{c}(n)=Y_1^{c}+\cdots+Y_n^{c}$ and $T^{c}(n)=J_1^{c}+\cdots+J_n^{c}$, we require that
$X^{c}(cu)\Rightarrow A(t)$ and $T^{c}(cu)\Rightarrow W(t)$ as $c\to\infty$, where the limits $A(t)$ and $W(t)$ are
 independent L\'evy processes.
Letting $N^{c}_t=\max\{n\geq 0:T^{c}(n)\leq t\}$, the CTRW scaling
limit $X^{c}(N_t^{c})\Rightarrow A(E^\nu_t)$
 see, \citet[Theorem 2.1]{M-S-triangular}.  A power-law mixture model for waiting times was proposed
  in \citet{M-S-ultra}:
   Take an i.i.d.\ sequence of mixing variables $(B_i)$ with $0<B_i<1$ and assume
   $P\{J_i^{c}>u|B_i=\beta\}=c^{-1}u^{-\beta}$ for $u\geq c^{-1/\beta}$, so that the waiting times are power
    laws conditional on the mixing variables.  The waiting time process $T^{c}(cu)\Rightarrow W(t)$ a
    nondecreasing L\'evy process, or subordinator, with $\E[e^{-s W(t)}]=e^{-t\psi_{W}(s)}$  and Laplace exponent
\begin{equation}\label{psiWdef}
\psi_W(s)=\int_0^\infty(e^{-s x}-1)\phi_{W}(dx) .
\end{equation}
The L\'evy measure
\begin{equation}\label{phiWdef}
\phi_W(t,\infty)=\int_0^1 t^{-\beta}\mu(d\beta)=\int_0^1 \frac{t^{-\beta}}{\Gamma(1-\beta)}\nu(d\beta),
\end{equation}
where $\mu$ is the distribution of the mixing variable, see \citet[Theorem 3.4 and Remark 5.1]{M-S-ultra}.  A
computation in  \citet[Eq.\ (3.18)]{M-S-ultra} using
$\int_0^\infty(1-e^{-s t})\beta
t^{-\beta-1}dt=\Gamma(1-\beta)s^\beta$ shows that
\begin{equation}\begin{split}\label{psiW}
\psi_{W}(s)
&= \int_0^1  s^\beta \Gamma(1-\beta) \mu(d\beta)= \int_0^1  s^\beta  \nu(d\beta) .
\end{split}\end{equation}
 Then $c^{-1}N^{c}_t\Rightarrow E^\nu(t)$, the inverse subordinator, see \citet[Theorem 3.10]{M-S-ultra}.  The general
 infinitely divisible L\'evy process limit $A(t)$ forms a strongly continuous convolution semigroup with
 generator $L$ (e.g., see \citet{ABHN}) and the corresponding CTRW scaling limit $A(E^\nu(t))$ is the stochastic
 solution to the distributed-order fractional Cauchy problem \citet[Eq.\
 (5.12)]{M-S-ultra} defined by
 \begin{equation}\label{DOFCPdef}
\D^{(\nu)} u(t,x)=L u(t,x).
\end{equation}
Since $\phi_W(0,\infty)=\infty$ in \eqref{phiWdef}, Theorem 3.1
in \citet{M-S-triangular} implies that the inverse subordinator
\begin{equation}\label{Etdef-D}
E^\nu(t)=\inf\{x>0:\ \ W(x)>t\}
\end{equation}
 has
density
\begin{equation}\label{etdensity-D}
g(t,x)=\int_0^t \phi_W(t-y,\infty) P_{W(x)}(dy) .
\end{equation}
This same condition  ensures also that $E^\nu(t)$ is almost surely
continuous, since $W(t)$ jumps in every interval, and hence is
strictly increasing. Further, it follows from the definition
\eqref{Etdef-D} that $E^\nu(t)$ is monotone nondecreasing.

We say that a function is a {\bf mild solution} to a pseudo-differential equation if its transform solves the corresponding equation in transform space. The next Lemma follows easily by taking Laplace transforms.

\begin{lemma}[\citet{mnv-2}]\label{eigenvalue-problem-D}
For any $\lambda>0$, $h(t, \lambda)=\int_0^\infty
e^{-\lambda x}g(t,x)\,dx=\E [e^{-\lambda E^\nu(t)}]$ is a  mild solution of
\begin{equation}\label{dist-order-density-pde-D}
\D^{(\nu)}h(t,\lambda)=-\lambda h(t, \lambda); \ \ h(0, \lambda)=1.
\end{equation}
\end{lemma}


\citet{koch4} considered the following: Let $\rho(\alpha)$ be a right continuous non-decreasing step function on $(0,1)$. Assume that $\rho$ has two sequences of jump points, $\beta_n$ and $\nu_n$, $n=0, 1,2,\cdots$, where $\beta_n\to 0,\  \nu_n\to 1$, $\beta_0=\nu_0\in (0,1)$. Suppose also that the sequence $\{\beta_n\}$ is strictly decreasing and $\{\nu_n\}$ is strictly increasing.
Let $\gamma^1_n=(\rho(\beta_n)-\rho(\beta_n-0))$ and $\gamma^2_n=(\rho(\nu_n)-\rho(\nu_n-0))$, and define the distributed order differential operator
 \begin{equation}\label{infinite-dist-def}
 \mathbb{D}^{(\rho)} u(t,x)=\sum_{n=0}^\infty
\gamma^1_n\partial ^{\beta_n}_tu(t,x)+\sum_{n=0}^\infty
\gamma^2_n\partial ^{\nu_n}_tu(t,x).
\end{equation}
Since $\rho$ is a finite measure we have
$$
 \sum_{n=1}^\infty\gamma^1_n<\infty, \ \ \sum_{n=1}^\infty\gamma^2_n<\infty.
$$

Here the corresponding subordinator is the sum of infinitely many independent stable subordinators;
\begin{equation}\label{infinite-dist-order-subordinator}
\begin{split}
W_t &=\sum_{n=0}^\infty
(\gamma^1_n)^{1/\beta_n}\Gamma(1-\beta_n)W_t^{\beta_n}\\
&+\sum_{n=1}^\infty
(\gamma^2_n)^{1/\nu_n}\Gamma(1-\nu_n)W_t^{\nu_n}
\end{split}
\end{equation}
for independent stable subordinators $W_t^{\beta_n}$, $W_t^{\nu_n}$ for $n=0,1,\cdots $.
\begin{lemma}\label{dist-rho-pde}Let $E^\rho(t)=\inf\{x>0:\ \ W(x)>t\}$. Then
$h(t,\mu)=\E (e^{-\mu E^\rho(t)})$ is the classical  solution to the eigenvalue problem
$$
\mathbb{D}^{(\rho)} h(t,\mu)=-\mu h(t,\mu), \ h(0,\mu)=1.
$$

\end{lemma}
Using inverse Laplace transforms  \cite{koch4} established the following representation of $h(t,\mu)$:
\begin{equation}\label{laplace-rep-hrho}
h(t,\mu)=\frac{\mu}{\pi}\int_0^\infty r^{-1}e^{-tr}\frac{H_1(r)}{H_2(r)}dr
\end{equation}
where
\begin{equation}\begin{split}
H_1(r)&=\sum_{n=0}^\infty\bigg[\gamma^1_n)r^{\beta_n}\sin(\pi\beta_n)+
(\gamma^2_n)r^{\nu_n}\sin(\pi\nu_n)\bigg]\\
H_2(r)&=\bigg\{\mu+\sum_{n=0}^\infty\bigg[\gamma^1_nr^{\beta_n}\cos(\pi\beta_n)+
\gamma^2_nr^{\nu_n}\cos(\pi\nu_n)\bigg]\bigg\}^2\\
&+\bigg\{\sum_{n=0}^\infty\bigg[\gamma^1_nr^{\beta_n}\sin(\pi\beta_n)+
\gamma^2_nr^{\nu_n}\sin(\pi\nu_n)\bigg]\bigg\}^2.
\end{split}\end{equation}

Let $D\subset \rd$ be a bounded domain with $\partial D \in
 C^{1,\alpha}$ for some $0<\alpha<1$, and  $D_\infty=(0,\infty )\times D$.
We will write $u\in C^k(\bar D)$ to mean that for each fixed
$t>0$, $u(t,\cdot)\in C^k(\bar D)$, and
 $u\in C_b^k(\bar D_\infty)$ to mean that $u\in C^k(\bar D_\infty)$ and is bounded.

\noindent Define
\begin{equation*}\begin{split}
\mathcal{H}_{\Delta}(D_\infty)&= \{u:D_\infty\to \rr :\ \ \Delta
u(t,x)\in C(D_\infty)\};\\
\mathcal{H}_{\Delta}^{b}(D_\infty)&= \mathcal{H}_{\Delta}(D_\infty)
\cap \{u: |\partial_t u(t,x)|\leq k(t)g(x),\ \\
&\ \ \ \ \ \  g\in
L^\infty(D),  \ t>0  \}, 
\end{split}\end{equation*}
for some functions  $k$ and $ b$
satisfying the condition
\begin{equation}\label{e-convolution-k}
b(\lambda)\int_0^1\int_0^t\frac{k(s)ds}{(t-s)^\beta}d\mu(\beta)<\infty,
\end{equation}
for $t, \lambda >0$ and
\begin{equation}\label{summability-k}
k(t)\sum_{n=1}^\infty b(\lambda_n)\bar{f}(n)|\phi_n(x)|<\infty.
\end{equation}
\begin{thm}
Let $f\in  C^1(\bar D)\cap C^2(D)$ for which the
  eigenfunction expansion (of $\Delta f$) with respect to the complete orthonormal basis $\{\phi_n:\ n\in \N \}$ converges uniformly and absolutely.
 Then the classical solution to the distributed-order fractional Cauchy problem
\begin{eqnarray}
\D^{(\rho)}u(t,x)
&=&
\Delta u(t,x),  \  \ x\in D, \ t\geq 0\label{infinite-e-frac-derivative-bounded-dd};\\
u(t,x)&=&0, \ x\in \partial D,\ t\geq 0; \nonumber\\
u(0,x)& =& f(x), \ x\in D,\nonumber
\end{eqnarray}
for  $u  \in
\mathcal{H}_{\Delta}^{b}(D_\infty)\cap C_b(\bar D_\infty) \cap C^1(\bar D)$, with the distributed order
fractional derivative $\D^{(\rho)}$ defined by \eqref{infinite-dist-def}, is given by
\begin{eqnarray}
u(t,x)&=&\E_{x}[f(B(E^\rho(t)))I( \tau_D(B)> E^\rho(t))]\nonumber\\
&=&\E_{x}[f(B(E^\rho({t})))I( \tau_D(B(E^\rho))> t)]\nonumber\\
&=& \int_{0}^{\infty}T_D(l)f(x)g(t,l)dl\nonumber\\
&=& \sum_0^\infty
\bar{f}(n)\phi_n(x)h(t, \mu_n).\label{infinite-e-stoch-rep-L}
\end{eqnarray}
In this case, $b(\lambda)=\lambda$, and  $k(t)$ is given by $k(t)=C t^{\beta_0-1}$,  $0<\beta_0<1$.
\end{thm}

\begin{pf}
The proof is similar to the proof of Theorem 4.1 in \citet{mnv-2}.
We give the main parts of the proof here.

Denote the Laplace transform $ t\rightarrow s $ of $u(t,x)$ by
\begin{equation*}\begin{split}
\tilde{u}(s,x)&=\int_{0}^{\infty}e^{-s t}u(t,x)dt .
\end{split}\end{equation*}
Since we are working on a bounded domain, the Fourier transform
methods in \cite{Zsolution} are not useful.  Instead, we will
employ Hilbert space methods.  Hence, given a complete orthonormal
basis $\{\phi_n(x)\}$ on $L^2(D)$, we will call
\begin{equation}\begin{split}
\bar{u}(t,n)&=\int_{D}\phi_n(x)u(t,x)dx;\\
\hat{u}(s,n)&=\int_{D}\phi_n(x)\int_{0}^{\infty}e^{-s t} u(t,x)dtdx   \\
&= \int_{D} \phi_{n}(x) \tilde{u}(s, x) dx   \\
& = \int_{0}^{\infty}e^{-s t} \bar{u}(t,x)dt~~\mbox{(when Fubini Thm holds)}
 \end{split}\end{equation}
respectively the $\phi_n$ and the $\phi_n$-Laplace transforms.
Since $\{\phi_n\}$ is a complete orthonormal basis for $L^2(D)$,
we can invert the $\phi_n$-transform to obtain
\[u(t,x)=\sum_n \bar{u}(t,n) \psi_n(x)\]
for any $t>0$, where the above series converges in the $L^2$ sense
(e.g., see \citet[Proposition 10.8.27]{Royden}).

Assume that $u(t,x)$ solves \eqref{infinite-e-stoch-rep-L}. Using Green's second identity, we obtain
$$
\int_D [u\Delta \phi_n-\phi_n\Delta u]dx=\int_{\partial
D}\left[u\frac{\partial \phi_n}{\partial \theta}-\phi_n\frac{\partial
u}{\partial \theta}\right]ds=0,
$$
since $u|_{\partial D}=0=\phi_n|_{\partial D}$,  $u\in
C^1(\bar{D})$  by assumption, and $\phi_n\in
C^1(\bar{D})$
  by \citet[Theorem 8.29]{gilbarg-trudinger}. Hence, the $\phi_n$-transform of $\Delta u$ is
\begin{equation}\begin{split}
\int_D \phi_n(x) \Delta u(t,x)  dx &
=
-\lambda_n\int_D  u(t,x) \phi_n(x) dx\\
& = -\lambda_n\bar{u}(t,n),
\end{split}\end{equation}
as $\phi_n$ is the eigenfunction of the Laplacian corresponding to
 eigenvalue $\lambda_n$.

 The fact that the operator $\D^{(\rho)}$ commutes with the $\phi_n$-transform follows from \eqref{e-convolution-k}.

Taking the $\phi_n$-transform of \eqref{infinite-e-stoch-rep-L} we obtain that
\begin{equation}\label{phi-trans-D}
\D^{(\rho)}\bar{u}(t,n)=
 -\lambda_n
\bar{u}(t,n).
\end{equation}

From Lemma \ref{dist-rho-pde} we get the solution
$$\bar u (t,n)=\bar f(n)h(t,\lambda_n)=\bar f(n)\E(e^{-\lambda_nE^\rho(t)}).$$

Now inverting the $\phi_n$-transform gives
$$u(t,x)= \sum_0^\infty
\bar{f}(n)\phi_n(x)h(t, \mu_n).
$$

The stochastic representation uses Lemma \ref{dist-rho-pde} and the stochastic representation of the killed semigroup of Brownian motion \eqref{stochastic-killed-bm}.

We use  the representation \eqref{laplace-rep-hrho} to establish
the fact that the solution is a classical solution. The details of the proof can be seen from the proof of the main results  in \citet{m-n-v-aop, mnv-2}.
\end{pf}

\begin{remark}\label{derivative-bound}

Suppose  that $\mu(d\beta)=p(\beta)d\beta$,
  the function $\beta\mapsto \Gamma(1-\beta)p(\beta)$ is in $C^1[0,1]$, $supp (\mu)=[\beta_0,\beta_1]\subset (0,1)$
   and $\mu(\beta_1)\neq 0$. Then
 \begin{equation}\label{kochubei-inverse}
h(t,\lambda)=\E [e^{-\lambda E^\nu(t)}]=\frac{\lambda}{\pi}\int_0^\infty r^{-1}e^{-tr}\frac{\Phi_1(r)}{\Phi_2(r)}dr
\end{equation}
where
\begin{equation*}
\begin{split}
\Phi_1(r)&=\int_0^1r^\beta\sin
(\beta\pi)\Gamma(1-\beta)p(\beta)d\beta\\
\Phi_2(r)&=[\int_0^1r^\beta \cos
(\beta\pi)\Gamma(1-\beta)p(\beta)d\beta+\lambda]^2\\
&+[\int_0^1r^\beta\sin
(\beta\pi) \Gamma(1-\beta)p(\beta)d\beta]^2.
\end{split}\end{equation*}

Suppose also that 
\begin{equation} \label{eqn3.2v3}
 C( \beta_0, \beta_1,p)= \int_{\beta_0}^{\beta_1}\sin
(\beta\pi)\Gamma(1-\beta)p(\beta)d\beta>0.
\end{equation}
Then
 $|\partial_t h(t, \lambda)|\leq \lambda k(t)$, where
\begin{equation}\label{time-bound-derivative}
k(t)=[C(\beta_0,\beta_1,p)\pi]^{-1}[\Gamma(1-\beta_1)t^{\beta_1-1}+\Gamma(1-\beta_0)t^{\beta_0-1}].
\end{equation}
In this case, $h(t,\lambda)$ is a classical solution to \eqref{dist-order-density-pde-D}.
The representation \eqref{kochubei-inverse} is due to  \citet{koch3}, which follows by inverting the
Laplace transform of  \eqref{dist-order-density-pde-D}.

For  $u  \in
\mathcal{H}^b_\Delta(D_\infty)\cap C_b(\bar D_\infty) \cap C^1(\bar D)$ for $k$ given by \eqref{time-bound-derivative},  \citet{mnv-2} shows that the solution  to \eqref{infinite-e-frac-derivative-bounded-dd}, where $\D^{(\rho)}$ replaced with the more general $\D^{(\nu)}$, is a strong (classical) solution  for $f\in C^1(\bar D)\cap C^2(D) $ for which $\Delta f$ has an absolutely and uniformly convergent eigenfunction expansion w.r.t $\{\phi_n\}$. \citet{naber} studied  distributed-order fractional Cauchy problem in $D=(0,M)\subset \mathbb R$.
\end{remark}
\begin{remark}\label{extension-L}
The methods of this paper also apply to the Cauchy problems that are obtained by replacing Laplacian with
 uniformly elliptic operator in divergence form  defined on
$C^2$ functions by
\begin{equation}\label{unif-elliptic-op}
Lu=\sum_{i,j=1}^{d}\frac{\partial \left(a_{ij}(x)(\partial
u/\partial x_i)\right)}{\partial x_j}
\end{equation}
with $a_{ij}(x)=a_{ji}(x)$ and, for some $\lambda>0,$
\begin{equation}\label{elliptic-bounds}
\lambda \sum_{i=1}^ny_i^2\leq \sum_{i,j=1}^na_{ij}(x)y_iy_j\leq
\lambda^{-1} \sum_{i=1}^ny_i^2,\ \ \forall y \in \rd.
\end{equation}
If $X_t$ is a solution to
$
dX_t=\sigma (X_t)dB_t+b(X_t)dt, \ \ X_0=x_0,
$
where $\sigma$ is a $d\times d$ matrix, and $B_t$ is a Brownian
motion, then $X_t$ is associated with the operator $L$ with
$a=\sigma \sigma ^T$, see Chapters 1 and 5 in
 \citet{bass}. Define the first exit time as
 $\tau_D(X)=\inf \{ t\geq 0:\ X_t\notin D\}$. The semigroup defined by
$T_D(t)f(x)=E_x[f(X_t)I(\tau_D(X))>t)]$ has generator $L$ with Dirichlet boundary conditions, which
follows by an application of the It$\mathrm{\hat{ o}}$ formula.
\end{remark}


\bibliography{erkan-ifacconf}             









%
\end{document}